\documentclass[10pt]{article}
\usepackage{amssymb}
\usepackage{latexsym}
\usepackage{amsmath,amsthm,xspace}

\newcommand{\streep}{\ |\ }
\newcommand{\Me}{\ensuremath{(M_{1},\del_{1})}\ }
\newcommand{\Mt}{\ensuremath{(M_{2},\del_{2})}\ }
\newcommand{\szwak}{\stackrel{\sigma w}{\longrightarrow}}
\newcommand{\zwak}{\stackrel{w}{\longrightarrow}}
\newcommand{\BH}{\ensuremath{\mbox{B}(H)}}
\newcommand{\BHL}{\ensuremath{\mbox{B}(H\otimes L)}}
\newcommand{\BHe}{\ensuremath{\mbox{B}(H_{1})}}
\newcommand{\BHt}{\ensuremath{\mbox{B}(H_{2})}}
\newcommand{\ten}{\otimes \iota}
\newcommand{\tenn}{\otimes \iota \otimes \iota}

\newcommand{\norm}[1]{\ensuremath{\Vert #1 \Vert}}
\newcommand{\C}{\ensuremath{\mathbb{C}}}
\newcommand{\R}{\ensuremath{\mathbb{R}}}

\newcommand{\cU}{\ensuremath{\mathcal{U}}}
\newcommand{\cV}{\ensuremath{\mathcal{V}}}

\newcommand{\cros}{M \, \mbox{$_{\alpha,\cU}$}\hspace{-.2ex}\mbox{$\ltimes$} \, N}
\newcommand{\cross}{M_{1} \, \mbox{$_{\alpha,\cU}$}\hspace{-.2ex}\mbox{$\ltimes$} \, M_{2}}
\newcommand{\kruisje}[1]{\, \mbox{$_{#1}$}\hspace{-.2ex}\mbox{$\ltimes$} \,}

\newcommand{\del}{\Delta}
\newcommand{\de}{\Delta}
\newcommand{\delh}{\hat{\Delta}}

\newcommand{\delopt}{\del_2 \hspace{-1ex}\raisebox{1ex}[0pt][0pt]{\scriptsize\fontshape{n}\selectfont op}}

\newcommand{\delop}{\del \hspace{-.3ex}\raisebox{0.9ex}[0pt][0pt]{\scriptsize\fontshape{n}\selectfont op}}
\newcommand{\delhop}{\delh \hspace{-.3ex}\raisebox{0.9ex}[0pt][0pt]{\scriptsize\fontshape{n}\selectfont op}}
\newcommand{\tekst}[1]{\;\text{#1}\;}
\newcommand{\nsf}{n.s.f.\xspace\!}
\newcommand{\vfi}{\varphi}

\newcommand{\Nfi}{{\cal N}_{\vfi}}
\newcommand{\strong}{\mbox{$\sigma$-strong$^*$}\xspace}
\newcommand{\sluit}{{- \; \strong}}
\newcommand{\Mhh}{\hat{M}
\hspace{-1.05ex}\hat{\rule{0ex}{2.0ex}}\hspace{1.05ex}}
\newcommand{\dehh}{\hat{\del}
\hspace{-.95ex}\hat{\rule{0ex}{2.05ex}}\hspace{.95ex}}

\newcommand{\halpha}{\hat{\alpha}}
\newcommand{\hbeta}{\hat{\beta}}
\newcommand{\Proof}{$Proof.$\ }
\newcommand{\Endp}{\begin{flushright} $\qed$ \end{flushright}}

{\theoremstyle{definition}
\newtheorem{defi}{Definition}
\newtheorem{nota}[defi]{Notation}

\theoremstyle{plain}
\newtheorem{theorem}[defi]{Theorem}
\newtheorem{prop}[defi]{Proposition}

\newtheorem{lemma}[defi]{Lemma}

\setlength{\parindent}{0pt} \setlength{\parskip}{5pt plus 2pt minus 1pt} \frenchspacing
\sloppy

\begin{document}
\begin{center}
\Large\bf Amenability and the bicrossed product construction
\end{center}

\bigskip

\begin{center}
{\bf Pieter Desmedt and Johan Quaegebeur}

\smallskip

Department of Mathematics, K.U.Leuven \\
Celestijnenlaan 200 B \\
B-3001 Leuven (Belgium)

\smallskip

\begin{tabular}{ll}
e-mail : &\hspace{-0.3cm}Pieter.Desmedt@wis.kuleuven.ac.be  \\
 &\hspace{-0.3cm}Johan.Quaegebeur@wis.kuleuven.ac.be
\end{tabular}

\bigskip

{\bf Stefaan Vaes}\footnote{Research Assistant of the
Fund for Scientific Research -- Flanders (Belgium)\ (F.W.O.)} \smallskip\\
Institut de Math{\'e}matiques de Jussieu \\ Alg{\`e}bres d'op{\'e}rateurs et
repr{\'e}sentations, Plateau 7E \\ 175, rue du Chevaleret \\ F-75013 Paris (France)

\smallskip

\begin{tabular}{ll}
e-mail : &\hspace{-0.3cm}vaes@math.jussieu.fr
\end{tabular}

\bigskip

November 2001

\end{center}

\begin{abstract} \noindent
We study stability properties of amenable locally compact quantum groups under the bicrossed product construction. We
obtain as our main result an equivalence between amenability of the bicrossed product and amenability of the matched
quantum groups used as building ingredients of the bicrossed product. Finally, we give examples of non-amenable locally
compact quantum groups obtained by a bicrossed product construction.
\end{abstract}
\section{Introduction}
The theory of locally compact quantum groups has been introduced by
J.~Kustermans and the third author in \cite{KV1,KV2}, unifying compact
quantum groups and Kac algebras. As the example of the quantum
$SU_q(2)$-group, developed by Woronowicz, shows, the antipode of a
compact quantum group need not be bounded and it need not respect the
$^*$-operation. For this reason, compact quantum groups are not always
Kac
algebras. The crucial difference between Kac algebras and locally
compact quantum groups is the possible unboundedness of the antipode.

Taking into account the importance of amenable locally compact groups
within the category of all locally compact groups, it is natural to
consider amenability of locally compact quantum groups. In fact, the
main results on amenability of Kac algebras, have been developed by
Enock and Schwartz \cite{ES}, and their proofs can be repeated in the
more general framework of locally compact quantum groups. However,
there is still one open problem. Recall that there are many different
characterizations of amenability of locally compact groups. A first characterizations deals with the existence of an invariant
mean on a suitable algebra of functions on the group $G$
($L^\infty(G)$, or bounded continuous functions). Another
characterization says that the trivial representation of $G$ is weakly
contained in the left regular representation. Other
characterizations are most of the time closely related to one of these
two definitions. These two properties can be formulated for locally
compact quantum groups and, in this way, one defines
\emph{amenable} and \emph{strongly amenable} locally compact quantum
groups. It is known that all strongly amenable quantum groups are
amenable, but the converse has only been proven for locally compact
groups, see e.g.\ \cite{Green}, and for discrete Kac algebras
\cite{Ruan}.

Having defined amenability of locally compact quantum groups, one asks
for examples. A systematic way of constructing examples of locally
compact quantum groups has been developed by Majid \cite{Majid}, Baaj
and Skandalis \cite{B-S2} and Vainerman and the third author
\cite{VV}, and in this paper, we precisely characterize when these
locally compact quantum groups are amenable. We also give two examples
of non-amenable locally compact quantum groups obtained by this
so-called bicrossed product construction.

In \cite{VV}, one also defines bicrossed products of quantum groups,
and one makes the link with short exact sequences of locally compact
quantum groups, called extensions. In this paper, we will characterize
in this full generality, when the bicrossed product is amenable, and
in fact, our result is a quantum version of the well known result that
a locally compact group $G$ with normal closed subgroup $H$ is
amenable if and only if $H$ and $G/H$ are amenable.

\section{Preliminaries}


We refer to \cite{KV1} and \cite{KV2} for the theory of locally compact quantum groups in the $C^{*}$-algebra, as well
as in the von Neumann algebra language. For the non-specialists,
\cite{thesis} is a good starting point. We recall from
\cite{KV2} the definition of a von Neumann algebraic quantum group: $(M,\del)$ is called a \emph{(von Neumann
algebraic) locally compact quantum group} when
\begin{itemize}
\item $M$ is a von Neumann algebra and $\del :M \to M\otimes M$ is a normal and unital
$*$-homomorphism satisfying the coassociativity relation : $(\del \ten)\del = (\iota \otimes
\del)\del$;
\item there exist normal, semi-finite, faithful (n.s.f.) weights $\varphi$ and $\psi$ on $M$
such that
\begin{itemize}
\item $\varphi$ is left invariant, i.e. $\varphi\bigl((\omega \ten)\del(x)\bigr)=\varphi(x)\omega(1)$
for all $x \in \mathcal{M}^{+}_{\varphi}$ and $\omega \in M^{+}_{*}$,
\item $\psi$ is right invariant, i.e. $\psi\bigl((\iota \otimes \omega)\del(x)\bigr)=\psi(x)\omega(1)$
for all $x \in \mathcal{M}^{+}_{\psi}$ and $\omega \in M^{+}_{*}$.
\end{itemize}
\end{itemize}
Here, we use the notation $\mathcal{M}^{+}_{\varphi} = \{x \in M^+
\mid \varphi(x) < +\infty \}$, and analogously for
$\mathcal{M}^{+}_{\psi}$.

Fix a left invariant \nsf weight $\varphi$ on $(M,\del)$ and represent $M$ on the GNS-space of
$\varphi$ such that $(H,\iota,\Lambda)$ is a GNS-construction for $\varphi$. Then, we can
define a unitary $W$ on $H \otimes H$ by $$W^* (\Lambda(a) \otimes \Lambda(b)) = (\Lambda
\otimes \Lambda)(\del(b)(a \otimes 1)) \quad\text{for all}\; a,b \in \Nfi \; .$$ Here,
$\Lambda \otimes \Lambda$ denotes the canonical GNS-map for the tensor product weight $\vfi
\otimes \vfi$. One proves that $W$ satisfies the pentagonal equation: $W_{12} W_{13} W_{23} =
W_{23} W_{12}$. We say that $W$ is a \emph{multiplicative unitary}. The comultiplication can be given
in terms of $W$ by the formula $\del(x) = W^* (1 \otimes x) W$ for all $x \in M$. Also the von
Neumann algebra $M$ can be written in terms of $W$ as $$M = \{ (\iota \otimes \omega)(W) \mid
\omega \in \BH_* \}^\sluit \; .$$

Next, the locally compact quantum group $(M,\del)$ has an \emph{antipode} $S$, which is the unique \strong closed
linear map from $M$ to $M$ satisfying \linebreak $(\iota \otimes \omega)(W) \in D(S)$ for all $\omega \in \BH_*$,
$S(\iota \otimes \omega)(W) = (\iota \otimes \omega)(W^*)$ and such that the elements $(\iota \otimes \omega)(W)$ form
a \strong core for $S$. The antipode $S$ has a polar decomposition $S = R \tau_{-i/2}$, where $R$ is an anti-automorphism of $M$ and
$(\tau_t)$ is a strongly continuous one-parameter group of automorphisms of $M$. We call $R$ the \emph{unitary
antipode} and $(\tau_t)$ the \emph{scaling group} of $(M,\del)$. It is known that $\sigma (R \otimes R) \del = \del R$,
where $\sigma$ denotes the flip map on $M \otimes M$.

We turn the predual $M_*$ into a Banach algebra with product $\omega \ast \mu = (\omega
\otimes \mu)\del$, for all $\omega,\mu \in M_*$.

We use the notation $\delop$ to denote the opposite comultiplication defined by $\delop:=
\sigma \del$.

The dual locally compact quantum group $(\hat{M},\delh)$ is defined as follows. Its von Neumann algebra $\hat{M}$ is
$$\hat{M} = \{(\omega \otimes \iota)(W) \mid \omega \in \BH_* \}^\sluit$$ and the comultiplication is given by
$\delh(x) = \Sigma W (x \otimes 1) W^* \Sigma$ for all $x \in \hat{M}$, where $\Sigma$ denotes the flip map on the
tensorproduct of Hilbert spaces.

Since $(\hat{M},\delh)$ is again a locally compact quantum group, we can introduce the
antipode $\hat{S}$, the unitary antipode $\hat{R}$ and the scaling group $(\hat{\tau}_t)$
exactly as we did it for $(M,\del)$. Also, we can again construct the dual of
$(\hat{M},\delh)$, starting from the left invariant weight $\hat{\vfi}$ with GNS-construction
$(H,\iota,\hat{\Lambda})$. From the biduality theorem, we get that the bidual locally compact
quantum group $(\Mhh,\dehh)$ is isomorphic to $(M,\del)$.

Define $M_{c}$ to be the norm closure of the space \[\{(\iota \otimes \omega)(W) \quad |\quad \omega \in \BH_{*}\}\]
and $\del_{c}$ to be the restriction of $\del$ to $M_{c}$. It is proven in \cite{KV2} that the pair $(M_{c},\del_{c})$
is a reduced $C^{*}$-algebraic locally compact quantum group. We know that there is a bijective correspondence between
reduced $C^{*}$-algebraic quantum groups and von Neumann algebraic quantum groups. So, the choice for the von Neumann
algebra language is not a restriction.

A $*$-homomorphism
$\varepsilon: M_c \to \C$ is called a \emph{co-unit} of $(M_c,\de_c)$, if
\[(\varepsilon \otimes \iota)\del=(\iota \otimes
\varepsilon)\del=\iota \; .\]

Classical locally compact groups appear as $M = L^\infty(G)$ with
$\de(f)(p,q) = f(pq)$. The invariant weights are defined by
integrating with respect to the left or the right Haar measure.
The dual $\hat{M}$ can be identified with the group von Neumann
algebra $\mathcal{L}(G)$.

Working with tensor products with more than two factors, we will sometimes use the leg-numbering notation. For example,
if $H, K$ and $L$ are Hilbert spaces and $X \in \BHL$, we denote by $X_{13}$ (respectively, $X_{12},\ X_{23}$) the
operator $(1 \otimes \Sigma^*)(X \otimes 1)(1 \otimes \Sigma)$ (respectively, $X\otimes 1,\ 1\otimes X$) defined on $H
\otimes K \otimes L$. If now $H = H_1 \otimes H_2$ is itself a tensor product of two Hilbert spaces, then we sometimes
switch from the leg-numbering notation with respect to $H \otimes K \otimes L$ to the one with respect to the finer
tensor product $H_1 \otimes H_2\otimes K \otimes L$, for example, from $X_{13}$ to $X_{124}$. There is no confusion
here, because the number of legs changes. Weak and $\sigma$-weak convergence are denoted by $\zwak$, respectively
$\szwak$.

\section{Amenability}
Let $(M,\del)$ be a von Neumann algebraic locally compact quantum
group. A state $m \in M^{*}$ is said to be a \emph{left invariant mean (LIM)} on $(M,\del)$ if
\[m((\omega \otimes \iota)\del(x))=m(x)\omega(1),\] for all $\omega \in M_{*}$ and $x \in M$.
It is said to be a \emph{right invariant mean (RIM)} if
\[m((\iota \otimes \omega )\del(x))=m(x)\omega(1),\] for all $\omega \in M_{*}$ and $x \in M$.
Finally, if $m$ is both a LIM and a RIM, we call $m$ an \emph{invariant mean (IM)}.
\begin{defi}
We call $(M,\del)$ \emph{amenable} if there exists a left invariant mean (LIM) on $(M,\del)$.
We say that $(M,\del)$ is \emph{coamenable} if $(\hat{M},\delh)$ is amenable.
\end{defi}

\begin{defi}
We call $(M,\del)$ \emph{strongly amenable} if there exists a bounded co-unit on $(\hat{M}_{c},\delh_{c})$.
\end{defi}

In the Preliminaries section, we saw that classical locally
compact groups appear as $L^\infty(G)$ in the theory of locally
compact quantum groups. We have defined amenability in such a way
that, for every locally compact group $G$, we have that $G$ is
amenable if and only if $(L^{\infty}(G),\del_{G})$ is amenable.
Other authors sometimes use a ``dual'' terminology. This
difference originates from the choice which quantum group is
associated with a locally compact group, $L^{\infty}(G)$ or
$\mathcal{L}(G)$. Here, we adopt the point of view of Enock and
Schwartz \cite{ES} and Ruan \cite{Ruan} (i.e. we take
$L^{\infty}(G)$ as the associated quantum group). The ``dual''
convention is used, amongst others, by Banica \cite{Ba}, Baaj and
Skandalis \cite{B-S2}, Ng \cite{Ng} and Bedos, Murphy and Tuset
\cite{Tuset}. They use \emph{coamenable} where we use
\emph{strongly amenable}. So, our notion of ``coamenability''
disagrees with their notion. Whenever we cite a result of one of
the papers mentioned with different terminology, it will be
already translated to our setting.

M. Enock and J.M. Schwartz prove in \cite{ES} that, for Kac algebras, the following statements are equivalent with the
fact that a Kac algebra is strongly amenable:
\begin{itemize}
\item[(i)] There exists a net $(\xi_{j})_{j}$ of normalized vectors in $H$ such that \[ (\iota
\otimes \omega_{\xi_{j}})(W) \zwak 1,\]
\item[(ii)] There exists a bounded left (resp., right) approximate unit on $\hat{M}_{*}$.
\end{itemize}
It is proven in \cite{ES} that strong amenability implies amenability. They also claim that the opposite implication is
true, but, as mentioned by Ruan \cite{Ruan}, there is a gap in their proof. It is an important open question whether or
not amenability implies strong amenability. Until now, this is only known to be true for locally compact groups, see
for example \cite{Green}, and for discrete Kac algebras \cite{Ruan2}.

Further they show that the following statements are equivalent:
\begin{itemize}
\item[(i)] there exist a LIM on $(M,\del)$ (resp., RIM);
\item[(ii)] there exists a net $(\omega_{i})_{i}$ of states in $M_{*}$ such that $\omega \ast
\omega_{i}-\omega_{i}$ converges weakly to $0$ (resp., $\omega_{i}
\ast \omega - \omega_{i}$), for all $\omega \in M_*$ with $\omega(1)=1$;
\item[(iii)] there exists a net $(\omega_{i})_{i}$ of states in $M_{*}$ such that $\norm{\omega \ast
\omega_{i}-\omega_{i}}$ converges to $0$ (resp., $\norm{\omega_{i}
\ast \omega - \omega_{i}}$), for all $\omega \in M_*$ with $\omega(1)=1$.
\end{itemize}

All these results are also true for locally compact quantum groups. Not surprisingly, we can prove the following
proposition.
\begin{prop}
Let $(M,\del)$ be a locally compact quantum group. There exists a LIM on $(M,\del)$ if and only if there exists an
invariant mean on $(M,\del)$.
\end{prop}
\Proof One implication is immediate.

Conversely, suppose there exists a LIM on $(M,\del)$. From the
result mentioned above, we know that there exists a net of states
$(\omega_{i})_{i}$ in $M_{*}$ such that $\norm{\omega
  \ast \omega_{i}-\omega_{i}}$ converges to $0$ for all $\omega \in
M_*$ with $\omega(1)=1$. It is obvious that this is equivalent
with the existence of a net of states $(\omega^{\circ}_{i})_{i}$
in $M_{*}$ such that $\norm{\omega^{\circ}_{i} \ast \omega
-\omega^{\circ}_{i}}$ converges to $0$ for all $\omega \in M_*$
with $\omega(1)=1$, take $\omega^{\circ}_{i}=\omega_{i}\circ R$.
It is easy to prove that $\mu_{k}=(\omega_{i} \ast
\omega^\circ_{j})_{(i,j)}$ is a net of states such that
\[\norm{\mu_{k}\ast \omega - \mu_{k}}\to 0\ \tekst{and}\ \norm{\omega \ast \mu_{k} - \mu_{k}}\to
0\] for all $\omega \in M_{*}$ with $\omega(1)=1$. Let $m$ be a weak-$*$ limit point of $(\mu_{k})_{k}$ in the unit ball of $M^{*}$. It is
obvious that $m$ will be an invariant mean.

\Endp

\section{Bicrossed products} In this section, we collect some results and definitions treated in \cite{VV}.
\begin{defi}
We call a pair $(\alpha,\cU)$ a \emph{cocycle action} of a locally compact quantum group
$(M,\del)$ on a von Neumann algebra $N$ if
\[\alpha: N \to M\otimes N\]
is a normal, injective and unital $*$-homomorphism,
\[\cU \in M\otimes M \otimes N\] is a unitary, and if $\alpha$ and $\cU$ satisfy
\begin{align*}
(\iota \otimes \alpha)\alpha(x) &= \cU\, (\del \otimes \iota)\alpha(x)\, \cU^{*}\, \tekst{for
all}\ x\in N,\\ (\iota \otimes \iota \otimes \alpha)(\cU)(\del \otimes \iota \otimes
\iota)(\cU) &= (1\otimes \cU)(\iota \otimes \del \otimes \iota)(\cU).
\end{align*}
\end{defi}
If $\cU$ is trivial, i.e. $\cU = 1$, we call $\alpha$ an \emph{action}.
\begin{nota} \label{12}
If $(\alpha,\cU)$ is a cocycle action of $(M,\del)$ on a von Neumann algebra $N$, we introduce
the notation $$\tilde{W} = (W \otimes 1)\ \cU^*$$ and then, $\tilde{W}$ is a unitary in $M
\otimes \BH \otimes N$.
\end{nota}
Given a cocycle action $(\alpha, \cU)$ of a locally compact quantum group $(M, \del)$ on a von Neumann algebra N, we
construct the crossed product $\cros$. This is the von Neumann subalgebra of $\BH \otimes N$ generated by \[\alpha(N)
\tekst{and}  \{(\omega \tenn)((W\otimes 1)\ \cU^{*}) \streep \omega \in M_{*}\}.\] When $\cU$ is trivial, the crossed
product is denoted by $M \kruisje{\alpha} N$. There is a unique action $\halpha$ of $(\hat{M},\delhop)$ on $\cros$ such
that, for all $x \in N$,
\begin{equation} \label{334}
\halpha\bigl(\alpha(x)\bigr)=1\otimes \alpha(x) \tekst{and} (\iota \otimes
\halpha)\bigl((W\otimes 1)\ \cU^{*}\bigr)=W_{12}W_{13}\cU^{*}_{134}.
\end{equation}
We call this action $\halpha$ the dual action. It is proven in \cite{VV} that the fixed point
algebra
\[(\cros)^{\halpha}=\{x \in \cros \streep \halpha(x)=1\otimes x \}=\alpha(N).\]

\begin{defi}
A pair $(M_{1},\del_{1}),(M_{2},\del_{2})$ is said to be a \emph{matched pair} of locally compact quantum groups if
there exists a triple $(\tau,\cU,\cV)$ (called a \emph{cocycle matching}) satisfying the following conditions:
\begin{itemize}
\item $\cU \in M_{1}\otimes M_{1} \otimes M_{2}$ and $\cV \in M_{1} \otimes M_{2} \otimes
M_{2}$ are both unitaries,
\item $\tau: M_{1}\otimes M_{2} \to M_{1}\otimes M_{2}$ is a faithful $*$-homomorphism,
\item defining $\alpha(y)=\tau(1\otimes y)$ and $\beta(x)=\tau(x\otimes
1)$ we have \begin{itemize}
\item $(\alpha,\cU)$ is a cocycle action of $(M_{1},\del_{1})$ on $M_{2}$,
\item $(\sigma\beta,\cV_{321})$ is a cocycle action of $(M_{2},\del_{2})$ on $M_{1}$,
\item $(\alpha,\cU)$ and $(\beta,\cV)$ are matched in the following sense:
\begin{align*}
\tau_{13}(\alpha \otimes \iota)\del_{2}(y) &= \cV_{132}(\iota \otimes \del_{2})\alpha(y)
\cV_{132}^{*},\\ \tau_{23}\sigma_{23}(\beta \otimes \iota)\del_{1}(x)&=\cU (\del_{1}\otimes
\iota)\beta(x)\ \cU^{*},\\ (\del_{1} \tenn)(\cV)(\iota \ten \otimes \delopt)(\cU^{*}) &= \\
(\cU^{*}\otimes 1)(\iota \otimes \tau \sigma \otimes \iota)\bigl(&(\beta \tenn)(\cU^{*})(\iota
\ten \otimes \alpha)(\cV)\bigr)(1 \otimes \cV).
\end{align*}
\end{itemize}
\end{itemize}
\end{defi}
Given a cocycle matching  $(\tau,\cU,\cV)$ of $(M_{1},\del_{1})$ and
$(M_{2},\del_{2})$, one is able to construct the
cocycle bicrossed product $(M,\del)$. By definition $M=\cross$ and $\del(x)=W^{*}(1\otimes x)W$ with $W=\Sigma
\hat{W}^{*} \Sigma$ and \begin{align*} \hat{W}&=(\beta \tenn)\bigl((W_{1}\otimes 1)\ \cU^{*}\bigr)(\iota \ten \otimes
\alpha)\bigl(\cV(1 \otimes \hat{W}_{2})\bigr)\\ &\in M_{1}\otimes \BHt
\otimes \BHe \otimes M_{2}. \end{align*} It is
proven in \cite{VV} that $(M,\de)$ is a locally compact quantum group
and that $W$ is its multiplicative unitary.

In Section~5, \Me and \Mt will always be two locally compact quantum groups matched by
$(\tau,\cU,\cV)$ and their cocycle bicrossed product locally compact quantum group will be denoted by
$(M,\del)$. All the objects associated with a quantum group (e.g. $W$, $\del$,\ldots) will be
denoted with an index, when they refer to \Me and \Mt respectively and without an index when
they refer to $(M,\del)$. So we have that $W_{1}$ (resp., $W_{2}$) is the multiplicative
unitary of \Me (resp., \Mt) and
\[\tilde{W}_{1}=(W_{1}\otimes 1)\ \cU^{*}.\]

From propositions 2.4 and 2.5 of \cite{VV}, we know how the comultiplication $\del$ works on the generators $\alpha(x)$
and $(\omega \otimes \iota \otimes \iota)(\tilde{W})$ of $(M,\del)$.
\begin{align}
\del(\alpha(x))&=(\alpha \otimes \alpha)\del_{2}(x) \notag\\ (\iota \otimes
\delop)(\tilde{W}_1)&=(\tilde{W}_1\otimes 1\otimes 1)((\iota \otimes \alpha)\beta \otimes \iota
\otimes \iota)(\tilde{W}_1)(\iota \otimes \alpha \otimes \alpha)(\cV)\label{355}
\end{align}
Define $\hat{M}$ as the von Neumann subalgebra of $M_{1}\otimes B(H_{2})$ generated by $\beta(M_{1})$ and $\{(\iota
\otimes \iota \otimes \omega)(\cV(1\otimes \hat{W}_{2})) \streep \omega \in M_{2\ast} \}$.  We
define $\delh(z)=\hat{W}^{\ast}(1\otimes z)\hat{W}$, for all $z \in \hat{M}$. It is proven in \cite{VV} that $(\hat{M},\delh)$ is the dual
locally compact quantum group of $(M,\del)$. So, if we interchange the roles of $\alpha$ and $\beta$, $M_{1}$ and $M_{2}$
respectively, then we find, as the cocycle bicrossed product, the dual of
the original cocycle bicrossed product.

\begin{defi} \label{defstab}
A cocycle action $(\alpha,\cU)$ of $(M,\del)$ on a von Neumann algebra $N$ is said to be \emph{stabilizable} with a
unitary $X\in M \otimes N$ if $$(1 \otimes X)(\iota \otimes \alpha)(X) = (\del \otimes \iota)(X) \cU^* \; .$$
\end{defi}
\begin{prop} \label{propstab}
Let $(\alpha,\cU)$ be a cocycle action of $(M,\del)$ on $N$ which is stabilizable with a unitary $X\in M \otimes N$.
Then the formulas $$\beta:N \to M \otimes N : \beta(x) = X \alpha(x) X^*\ \tekst{and}\ \Phi: z \mapsto X^* z X$$
define, respectively, an action of $(M,\del)$ on $N$ and a $^*$-isomorphism from $M \kruisje{\beta} N$ onto $\cros$
satisfying
$$\halpha \circ \Phi = (\iota \otimes \Phi) \circ \hat{\beta} \; .$$
\end{prop}
The next proposition shows that many cocycle actions are stabilizable.
\begin{prop} \label{propstab2}
Let $(\alpha,\cU)$ be a cocycle action of $(M,\del)$ on $N$. Then  \linebreak $(\alpha \otimes \iota, \cU \otimes 1)$
is a cocycle action of $(M,\del)$ on $N \otimes \BH$ which is stabilizable.
\end{prop}

\section{Amenability and the bicrossed product construction}
In this section, we investigate the relation between amenability of
the cocycle bicrossed product quantum group
$(M,\del)$ and of its building ingredients $(M_{1},\del_{1})$ and $(M_{2},\del_{2})$.

We start with a technical remark about slicing with non-normal functionals. Let $N$ and $L$ be
von Neumann algebras, $n\in N^{*}$ and $X\in N\otimes L$.

If $n \in N_{*}$, then it is obvious that $(n \otimes \iota)(X) \in L$. This remains
true for $n\in N^{*}$, even if $n$ is not normal. Indeed, consider the map $L_{*}\to\C: \omega
\mapsto n((\iota\otimes \omega)(X))$. It is obvious that this is a bounded linear functional
and since $L=(L_{*})^{*}$, we know that there exists a unique $Y\in L$ such that
$\omega(Y)=n((\iota\otimes \omega)(X)$ for all $\omega \in L_{*}$. Denote $Y=(n\otimes
\iota)(X)$.

Suppose that $\Phi:L\to K$ is a normal $\ast$-homomorphism of von Neumann algebras. Since for
all $\omega \in K_{*}$
\begin{align*}
\omega(\Phi(n\otimes \iota)(X)) &= n((\iota \otimes \omega \circ \Phi)(X)),\\
 &= n((\iota \otimes \omega)(\iota \otimes \Phi)(X)),\\
 &= \omega((n\otimes \iota)(\iota \otimes\Phi)(X)),
\end{align*}
we may conclude that for all $\Phi$
\[\Phi((n\otimes \iota)(X))=(n\otimes \iota)(\iota \otimes \Phi)(X).\]
This will be used several times in the sequel, where $n$ will be an invariant mean and
$\Phi=\alpha,\del,\ldots$.

\begin{defi}
If $\alpha$ is an action of $(M,\del)$ on a von Neumann algebra $N$, we define an \emph{$\alpha$-invariant mean} to be a
state $m \in N^{*}$ such that
\[m((\omega \otimes \iota)\alpha(x))=m(x)\omega(1)\] for all $\omega \in M_{*}$ and $x \in N$.
\end{defi}

\begin{prop} \label{gekruist}
Let $(\alpha, \mathcal{U})$ be a cocycle action of $(M,\del)$ on a von Neumann algebra $N$,\;
$\cros$ the cocycle crossed product and $\halpha$ the dual action. Then, $(\hat{M},\delh)$ is amenable
if and only if there exists a $\halpha$-invariant mean on $\cros$.
\end{prop}

\Proof Suppose that $\hat{m}$ is an invariant mean on $(\hat{M},\delh)$. Then we argue that there exists an
$\halpha$-invariant mean on $\cros$. This can be done by generalizing a result in \cite{ES} from the Kac algebra level
to the general setting. However, their proof is based on a non-constructive argument. We construct explicitly an
$\halpha$-invariant mean on $\cros$. The dual weight construction is the source of inspiration. Define $T:\cros \to
\cros$ by $T(z):= (\hat{m}\otimes \iota)\halpha(z)$. We prove that $T(z)\in \alpha(N)$ for all $z \in \cros$. Since
$\alpha(N)$ is the fixed point algebra of $\halpha$, it is sufficient to show that $\halpha(T(z))=1\otimes T(z)$.\\
Observe that, $\halpha$ being an action, $\halpha(T(z))=(\hat{m}\otimes \iota\otimes \iota)(\delhop\otimes \iota)\halpha(z)$.\\
So we have to prove that, for all $\omega \in (\hat{M}\otimes \cros)_{*}$,
\[\omega((\hat{m}\otimes \iota\otimes \iota)(\delhop\otimes
\iota)\halpha(z))=\omega(1\otimes(\hat{m}\otimes \iota)\halpha(z)).\] But, it is sufficient to
check this for normal functionals of the form $\mu\otimes \nu$ with $\mu \in \hat{M}_{*}$ and
$\nu \in (\cros)_{*}$. Using the fact that $\hat{m}$ is a LIM on $(\hat{M},\delh)$ and hence a
RIM on $(\hat{M},\delhop)$, we get
\begin{align*}
(\mu \otimes \nu)((\hat{m}\otimes \iota\otimes \iota)(\delhop\otimes \iota)\halpha(z)) &=
\hat{m}(( \iota\otimes \mu \otimes \nu)(\delhop\otimes \iota)\halpha(z)),\\ &=
\mu(1)\hat{m}((\iota\otimes \nu)(\halpha(z))),\\ &= (\mu \otimes \nu)(1\otimes(\hat{m}\otimes
\iota)\halpha(z)).
\end{align*}
So we may conclude that $T(z) \in \alpha(N)$ for all $z\in \cros$.

Choose a state $\eta \in N^{*}$. Define $m(z)=\eta(\alpha^{-1}(T(z)))$. We will prove that $m$
is $\halpha$-invariant. For all $\omega \in \hat{M}_{*}$ and $z \in \cros$, we get that
\begin{align*}
m((\omega \otimes \iota)\halpha(z)) &= \eta(\alpha^{-1}((\hat{m}\otimes \iota)\halpha((\omega
\otimes \iota)\halpha(z)))),\\ &= \eta(\alpha^{-1}((\hat{m}\otimes \iota)(\omega \otimes
\iota \otimes \iota)((\iota \otimes \halpha)\halpha (z)))),\\ &=
\eta(\alpha^{-1}((\hat{m}\otimes \iota)(\omega \otimes \iota \otimes \iota)((\delhop\otimes
\iota)\halpha (z)))),\\ &= \eta(\alpha^{-1}((\hat{m}\otimes
\iota)\halpha(z)))\omega(1)=m(z)\omega(1).
\end{align*}

Conversely, suppose that $m$ is a $\halpha$-invariant mean on
$\cros$. We have to prove that
$(\hat{M},\delh)$ is amenable. The proof is cut into three cases.

{\sc Case 1: }\emph{$\cU$ is trivial.}\\ We know that $M \kruisje{\alpha} N$ is generated by
$\alpha(N)$ and $\hat{M}\otimes \C$. Define $\hat{m}(\hat{x}):=m(\hat{x}\otimes 1)$ for all
$\hat{x} \in \hat{M}$. Using the formula $\halpha(\hat{x}\otimes 1)=\delhop(\hat{x})\otimes
1$, we get that for all $\omega \in \hat{M}_{*}$
\begin{align*}
\hat{m}((\omega \otimes \iota)\delhop(\hat{x})) &= m((\omega \otimes \iota \otimes \iota
)(\delhop(\hat{x})\otimes 1)),\\ &= m((\omega \otimes \iota \otimes \iota
)\halpha(\hat{x}\otimes 1)),\\ &= m(\hat{x}\otimes 1)\omega(1)=\hat{m}(\hat{x})\omega(1).
\end{align*}
So, we may conclude that $\hat{m}$ is a left invariant mean on $(\hat{M},\delhop)$.

{\sc Case 2: }\emph{$(\alpha,\cU)$ is stabilizable.}\\ We know from Proposition~\ref{propstab} that, in this case,
there exists an action $\beta$ of $(M,\del)$ on $N$ and a $*$-isomorphism
\[\Phi: M \kruisje{\beta} N \to \cros,\] such that $\halpha \circ \Phi=(\iota \otimes \Phi)\circ
\hbeta$.\\ Define $\tilde{m}:=m\circ \Phi$, then it is easy to prove that $\tilde{m}$ is $\hbeta$-invariant. Using
the first case, we may conclude that the restriction of $\tilde{m}$ to $\hat{M}$ will be a LIM on $(\hat{M},\delhop)$.

{\sc General case: }\emph{Arbitrary $(\alpha,\mathcal{U})$.}\\ In general, $(\alpha \otimes \iota,\mathcal{U} \otimes
1)$  will be a cocycle action of $(M,\del)$ on $N\otimes \BH$ and we know from Proposition~\ref{propstab2} that it will
be stabilizable. It is not too difficult to show that its corresponding cocycle crossed product factorizes as $(\cros)
\otimes \BH$, as well as the dual action, which is given by $\halpha \otimes \iota$.\\ Choose a normalized vector $\xi
\in H$. Then, we have for all $z \in (\cros) \otimes \BH$ that
\begin{align*}
(\iota \otimes m \otimes \omega_{\xi})((\halpha \otimes \iota)(z)) &= (\iota \otimes
m)\halpha((\iota \otimes \omega_{\xi})(z)),\\ &= m((\iota \otimes \omega_{\xi})(z))1,\\ &=
(m\otimes \omega_{\xi})(z)1.
\end{align*}
We find that $m\otimes \omega_{\xi}$ is $(\halpha \otimes \iota)$-invariant and from the
second case, we may conclude that $(\hat{M},\delh)$ is amenable.

\Endp

With this theorem in mind, we are going to prove our main result, generalizing a result of Ng
\cite{Ngex}. Ng proves in \cite{Ngex} that the bicrossed product with trivial cocycles of two locally
compact groups $G_{1}$ and $G_{2}$, is amenable if $G_{2}$ is amenable. Notice that for any
group $G_{1}$, $(L^{\infty}(G_{1}),\del_{1})$ is always coamenable.

To
prove our main theorem, we need a lemma. We can get this result from propositions 3.1 and 3.4 of \cite{VV}, but we have chosen to give
a straightforward proof.

\begin{lemma} \label{lemstab1}
Let $(\tau,\mathcal{U},\mathcal{V})$ be a cocycle matching of $(M_{1},\del_{1})$ and
$(M_{2},\del_{2})$ and let $(M,\del)$ be the cocycle bicrossed product. Then
\begin{equation}
(\iota\otimes \delop)\halpha(z)=(\halpha\otimes \iota)\delop(z) \ \mbox{for all}\ z \in M.
\end{equation}
\end{lemma}
\Proof It suffices to check it on the generators. Choose $x\in M_{2}$. Observe that
\begin{align*}
(\iota\otimes \delop)\halpha\bigl(\alpha(x)\bigr) &= (\iota\otimes \delop)\bigl(1\otimes
\alpha(x)\bigr),\\
 &= 1\otimes \delop\bigl(\alpha(x)\bigr),\\
 &= 1\otimes (\alpha \otimes \alpha)\delopt(x),\\
 &= (\halpha\otimes \iota)\bigl((\alpha \otimes \alpha)\delopt(x)\bigr),\\
 &= (\halpha\otimes \iota)\delop\bigl(\alpha(x)\bigr).
\end{align*}
Now, we will prove that $(\iota\otimes \iota \otimes \delop)(\iota\otimes
\halpha)(\tilde{W_{1}})=(\iota\otimes \halpha\otimes \iota)\bigl((\iota\otimes
\delop)(\tilde{W_{1}})\bigr)$.\ Using Equation~(\ref{334}) we get
\[(\iota\otimes \iota \otimes \delop)(\iota\otimes
\halpha)(\tilde{W_{1}})=(\iota\otimes \iota \otimes
\delop)\bigl((W_{1})_{12}(\tilde{W_{1}})_{134}\bigr).\] Finally, observe that, as operators on
$H_{1}\otimes H_{1}\otimes H_{1}\otimes H_{2} \otimes H_{1} \otimes H_{2}$,
\begin{align*}
\lefteqn{(\iota\otimes \halpha\otimes \iota)\bigl((\iota\otimes
\delop)(\tilde{W_{1}})\bigr)}\\ &= (\iota \otimes\halpha \otimes
\iota)\bigl((\tilde{W_{1}}\otimes 1\otimes 1)\bigl((\iota \otimes \alpha)\beta\otimes
\iota\otimes \iota \bigr)(\tilde{W_{1}})(\iota\otimes \alpha \otimes
\alpha)(\mathcal{V})\bigr),\\
 &= (W_{1})_{12}(\tilde{W_{1}})_{134}((\iota \otimes
\alpha)\beta\otimes \iota\otimes \iota)(\tilde{W_{1}})_{13456}(\iota\otimes \alpha \otimes
\alpha)(\mathcal{V})_{13456},\\
 &= (W_{1})_{12}(\iota\otimes \delop)(\tilde{W}_{1})_{13456},\\
 &= (W_{1})_{12}(\iota\otimes \iota \otimes \delop)\bigl((\tilde{W}_{1})_{134}\bigr),\\
 &= (\iota\otimes \iota \otimes
\delop)\bigl((W_{1})_{12}(\tilde{W_{1}})_{134}\bigr),
\end{align*}
where we used Equation~(\ref{355}) in the first line.

\Endp

\begin{theorem} \label{zwak}
Let $(\mathcal{\tau},\mathcal{U},\mathcal{V})$ be a cocycle matching of $(M_{1},\del_{1})$ and
$(M_{2},\del_{2})$ and let $(M,\del)$ be the cocycle bicrossed product. Then, $(M,\del)$ is
amenable if and only if $(\hat{M}_{1},\delh_{1})$ and $(M_{2},\del_{2})$ are amenable.
\end{theorem}


\Proof The proof is divided into three parts.

1) \emph{If $(M,\del)$ is amenable, then $(\hat{M}_{1},\delh_{1})$ is amenable.}\\ Let $m$ be
an invariant mean on $(M,\del)$. From Proposition~\ref{gekruist}, we know that it is sufficient
to show that $m$ is $\halpha$-invariant. If we apply $\iota \otimes \iota \otimes m$ on the
result in Lemma~\ref{lemstab1}, we get that, for all $z \in M$,
\[(\iota \otimes m)\halpha(z)\otimes 1=m(z)1\otimes 1. \]
From this we conclude that $(\iota \otimes m)\halpha(z)=m(z)1$ and therefore that $m$ is
$\halpha$-invariant.

2) \emph{If $(M,\del)$ is amenable, then $(M_{2},\del_{2})$ is amenable.}\\ Suppose that $m$ is a LIM on $(M,\del)$.
Define $m_{2} \in M_{2}^{*}$ by $m_{2}(x)=m(\alpha(x))$. Since $(\alpha\otimes \alpha)\del_{2}=\del\circ \alpha$ and
$M_{2*}=\{\omega\circ \alpha \streep \omega\in M_{*}\}$ it is obvious that $m_{2}$ will be a LIM on $(M_{2},\del_{2})$.

3) \emph{If $(\hat{M}_{1},\delh_{1})$ and $(M_{2},\del_{2})$ are amenable, then $(M,\del)$ is
amenable.}\\
Suppose that $\hat{m}_{1}$ and $m_{2}$ are left invariant means on $(\hat{M}_{1},\delh_{1})$
and $(M_{2},\del_{2})$ respectively.\\ Consider the dual action $\hat{\alpha}: M \to
\hat{M}_{1}\otimes M$. Define $T(z):= (\hat{m}_{1}\otimes \iota)\hat{\alpha}(z)$ for all $z
\in M$. From the proof of Proposition~\ref{gekruist}, we know that $T(z)\in \alpha(M_{2})$.

Define $m:=m_{2}\circ \alpha^{-1} \circ T$. We prove that $m$ is a left invariant mean on M.
Choose any $z \in M$. Applying $\hat{m}_{1}\otimes \iota \otimes \iota$ on both sides of the
result of Lemma~\ref{lemstab1} we get
\[\delop(T(z))=(T\otimes \iota)\delop(z).\]
So, we have for all $\mu \in M_{*}$ \[(\iota \otimes \mu)\delop(T(z))=T((\iota \otimes
\mu)\delop(z)).\] Take $y\in M_{2}$ such that $T(z)=\alpha(y)$.\\ Since $(\iota \otimes
\mu)\delop(\alpha(y))=(\iota \otimes \mu)(\alpha\otimes \alpha)\delopt(y)=\alpha((\iota
\otimes \mu \circ \alpha)\delopt(y))$ we get
\begin{equation} \label{T1}
T((\iota \otimes \mu)\delop(z))=\alpha((\iota \otimes \mu \circ \alpha)\delopt(y)).
\end{equation}
When we apply $m_{2}\circ \alpha^{-1}$ on both sides of Equation~(\ref{T1}), we get \[m((\iota \otimes
\mu)\delop(z))=m_{2}((\iota \otimes \mu \circ \alpha)\delopt(y)).\] Now we can use left invariance of $m_{2}$ and we
find for all $\mu \in M_{*}$
\begin{align*}
m((\iota \otimes \mu)\delop(z)) &= m_{2}((\iota \otimes \mu \circ \alpha)\delopt(y))\\
 &= \mu(\alpha(1))m_{2}(y)\\
 &= \mu(1)m_{2}(\alpha^{-1}(T(z)))\\
 &= \mu(1)m(z).
\end{align*}
Therefore, $m$ is a left invariant mean on $(M,\del)$.

\Endp

A natural question is whether or not the strong version of Theorem~\ref{zwak} is true, i.e., Theorem~\ref{zwak} with
amenability replaced by strong amenability. We can only give a partial
answer. First of all, it is not too difficult to see that
$(\hat{M}_{1},\delh_{1})$ is strongly amenable if $(M,\del)$ is. Just suppose that the net $(\hat{\mu}_{k})_{k}$ is an
approximate unit of $\hat{M}_{*}$. Define $\mu_{1k}:=\hat{\mu}_{k}\circ \beta$, then $(\mu_{1k})_{k}$ is an approximate
unit of $M_{1 \ast}$. So we arrive at the following proposition.

\begin{prop} \label{st}
Let $(\mathcal{\tau},\mathcal{U},\mathcal{V})$ be a cocycle matching of $(M_{1},\del_{1})$ and
$(M_{2},\del_{2})$ and let $(M,\del)$ be the cocycle bicrossed product. If $(M,\del)$ is
strongly amenable, then $(\hat{M}_{1},\delh_{1})$ is strongly amenable.
\end{prop}

Next, we can prove the strong version of Theorem~\ref{zwak} in the
case where the cocycles are trivial: $\cU=\cV=1$. We do not know
whether or not the same result holds with non-trivial cocycles.

\begin{theorem}
Let $(M,\del)$ be the bicrossed product of $(M_{1},\del_{1})$ and $(M_{2},\del_{2})$ with
trivial cocycles. Then, $(M,\del)$ is strongly amenable if and only if
$(\hat{M}_{1},\delh_{1})$ and $(M_{2},\del_{2})$ are strongly amenable.
\end{theorem}
\Proof We will first prove that if $(\hat{M}_{1},\delh_{1})$ and $(M_{2},\del_{2})$ are
strongly amenable, then $(M,\del)$ is strongly amenable.\\ Suppose that $(\omega_{i})_{i}$ is
a bounded two-sided approximate unit for $M_{1*}$. It is sufficient to show that
\begin{equation} \label{c} ((\omega_{i}\otimes \iota)\beta \otimes
\iota)(W_{1}) \szwak 1.
\end{equation}
Indeed, by definition, $\hat{W}=((\beta\otimes \iota)(W_{1})\otimes 1)(1\otimes (\iota \otimes
\alpha)(\hat{W}_{2}))\in M_{1}\otimes \BHt \otimes \BHe \otimes M_{2}$ and so
\[(\omega_{i} \otimes \iota \otimes \iota \otimes \iota)(\hat{W})=(((\omega_{i} \otimes
\iota)\beta \otimes \iota)(W_{1}) \otimes 1)((\iota \otimes \alpha)(\hat{W}_{2})).\] Using Equation~(\ref{c}), we get
that
\[(\omega_{i} \otimes \iota \otimes \iota \otimes \iota)(\hat{W})\szwak
(\iota\otimes \alpha)(\hat{W}_{2}).\] Because $(M_2,\de_2)$ is
strongly amenable, we can
take a net $(\xi_{j})_{j}$ of normalized vectors in
$H_{2}$ such that $(\mu_{\xi_{j}} \otimes \iota)(\hat{W}_{2})\szwak 1$.

Choose $\mu \in M_{*}$. Observe that for all $i, j$
\[|\mu((\mu_{\xi_{j}} \otimes \iota \otimes \iota)(\omega_{i} \otimes \iota \otimes \iota \otimes
\iota)(\hat{W}))|\leq \norm{(\iota\otimes \iota \otimes \mu)(\hat{W})}.\] Taking first the limit over $i$ and then over
$j$ we get
\[|\mu(1)|\leq \norm{(\iota\otimes \iota \otimes \mu)(\hat{W})}\]
Define $\hat{\varepsilon}((\iota \otimes \mu)(\hat{W}))=\mu(1)$. Thus, $\hat{\varepsilon}$ is
a bounded co-unit for $(\hat{M}_{c},\delh_{c})$.

It remains to prove (\ref{c}). The definition of matched pairs implies that
\[\tau_{23}\sigma_{23}(\beta\otimes \iota)\del_{1}(x)=(\del_{1}\otimes \iota)\beta(x).\]
When we apply $\omega_{i}\otimes \iota \otimes \iota$ on both sides we get
\begin{equation} \label{a}
\tau \sigma(((\omega_{i}\otimes \iota)\beta \otimes \iota)\del_{1}(x))=(\omega_{i}\otimes
\iota \otimes \iota)((\del_{1}\otimes \iota)\beta(x)).
\end{equation}
For all $\omega\in M_{1*}$ and $\nu \in M_{2*}$, we have that \[(\omega \otimes \nu)(\omega_{i}\otimes \iota \otimes
\iota)((\del_{1}\otimes \iota)\beta(x))=(\omega_{i}\ast \omega \otimes \nu)\beta(x) \to (\omega \otimes
\nu)(\beta(x)).\] By linearity and the fact that $(\omega_{i}\otimes \iota \otimes \iota)((\del_{1}\otimes
\iota)\beta(x))$ is uniformly bounded in $i$, we get that $(\omega_{i}\otimes \iota \otimes \iota)((\del_{1}\otimes
\iota)\beta(x)) \szwak \beta(x)$.\\ Using Equation~(\ref{a}) and the normality of $\tau\sigma$ we find that
\[\tau \sigma(((\omega_{i}\otimes \iota)\beta \otimes \iota)\del_{1}(x)) \szwak
\beta(x)=\tau \sigma (1\otimes x).\] Now, $\tau \sigma$ is an injective and normal
$\ast$-homomorphism and therefore it will be homeomorphic onto his image for the $\sigma$-weak
topology (\cite{Dix}, p. 60). From this, we get
\begin{equation} \label{b} ((\omega_{i}\otimes \iota)\beta
\otimes \iota)\del_{1}(x) \szwak 1\otimes x.
\end{equation}
When we apply $\beta \otimes \iota \otimes \iota$ on $(\del_{1}\otimes
\iota)(W_{1})=W_{1,13}W_{1,23}$ we get
\[((\beta \otimes \iota)\del_{1}\otimes
\iota)(W_{1})=((\beta \otimes \iota)(W_{1}))_{124}W_{1,34}\] and
\[(((\omega_{i}\otimes \iota)\beta \otimes \iota)\del_{1}\otimes
\iota)(W_{1})=(((\omega_{i}\otimes \iota)\beta \otimes \iota)(W_{1}))_{13}W_{1,23}.\] Using
Equation~(\ref{b}), we may conclude that
\[(((\omega_{i}\otimes \iota)\beta \otimes \iota)(W_{1}))_{13}W_{1,23} \szwak1\otimes W_{1}=W_{1,23}.\]
As $W_{1}$ is invertible, this implies that
\begin{equation}
((\omega_{i}\otimes \iota)\beta \otimes \iota)W_{1} \szwak 1.
\end{equation}
This concludes the first part of the proof.
\bigskip

By taking trivial cocycles in Proposition~\ref{st}, it is immediately clear that $(\hat{M}_{1},\delh_{1})$ is strongly
amenable, if $(M,\del)$ is strongly amenable.

It remains to show that if $(M,\del)$ is strongly amenable, then $(M_{2},\del_{2})$ is
strongly amenable. Using the biduality theorem, it is sufficient to prove that if
$(\hat{M},\delh)$ is strongly amenable, then $(M_{1},\del_{1})$ is strongly amenable. Suppose
that $(\omega_{i})_{i}$ is a bounded two-sided approximate unit for $M_{*}$. We know that now
\[M=(\alpha(M_{2})\cup \{(\omega \otimes \iota)(W_{1}) \otimes 1 \streep \omega \in M_{1*}\})''
\; .\] Using Equation~(\ref{355}), we get
\[(\iota \otimes \delop)(W_{1}\otimes 1)=(W_{1}\otimes 1 \otimes 1 \otimes 1)((\iota \otimes
\alpha)\beta \otimes \iota \otimes \iota)(W_{1} \otimes 1),\] so
\begin{equation}\label{j}
(\iota \otimes \iota \otimes \omega_{i})(\iota \otimes \delop)(W_{1}\otimes 1)=(W_{1}\otimes
1)(\iota \otimes \alpha)\beta((\iota \otimes \omega_{i})(W_{1}\otimes 1)).
\end{equation}
Using the fact that $(\omega_{i})_{i}$ is an approximate unit of $M_{\ast}$, we have
\[(\iota \otimes \iota \otimes \omega_{i})(\iota \otimes \delop)(W_{1}\otimes
1) \szwak W_{1}\otimes 1\] and thus, by Equation~(\ref{j})
\[(\iota \otimes \alpha)\beta((\iota \otimes \omega_{i})(W_{1}\otimes
1)) \szwak 1.\] But $(\iota \otimes \alpha)\beta$ is a normal and injective
$\ast$-homomorphism and therefore
\begin{equation} \label{k}
(\iota \otimes \omega_{i})(W_{1}\otimes 1)  \szwak 1.
\end{equation}
Define $\mu_{i} \in \hat{M}_{1*}$ such that $\mu_{i}(z)=\omega_{i}(z\otimes 1)$ for all $z \in
\hat{M}_{1}$, so $(\iota\otimes \mu_{i})(W_{1})=(\iota \otimes \omega_{i})(W_{1} \otimes 1)$.
Using Equation~(\ref{k}) we get that
\[(\iota\otimes \mu_{i})(W_{1})\szwak 1\] and this concludes the proof.

\Endp

\section{Examples}
In order to construct these examples, we rely on the extension procedure of locally compact quantum groups as developed
in \cite{B-S2,Majid,VV}. All the bicrossed product locally compact quantum groups in \cite{VV} are
amenable. That is easily seen, since the groups from which one starts in the examples are both amenable. We give
two examples of non-amenable locally compact quantum groups, obtained by a bicrossed product construction. From theorem
\ref{zwak}, we know that, if we take, as one of the ingredients, a non-amenable group, the bicrossed product locally
compact quantum group will be not amenable. In the first we take $SL_{2}(\R)$ as the non-amenable group and in the
second (a double cover of) $SU(1,1)$. It is a known that these groups are not amenable, since these are non-compact,
almost connected, semi-simple Lie groups, see \cite{Green}.

We briefly review what is needed from the extension procedure.

Let $G$, $G_{1}$ and $G_{2}$ be locally compact groups with fixed left invariant Haar
measures. Let $i:G_{1} \to G$ be a homomorphism and $j: G_{2} \to G$ an antihomomorphism such
that both have a closed image and are homeomorphisms onto these images. Suppose moreover that
the mapping
\[\theta:G_{1}\times G_{2}\to \Omega \subset G:(g,s)\mapsto i(g)j(s)\] is a homeomorphism of
$G_{1} \times G_{2}$ onto an open subset $\Omega$ of G having a complement of measure zero.
Then we call $G_{1}$ and $G_{2}$ a \emph{matched pair pair of locally compact groups}. From
this data, one constructs a cocycle matching of $(L^{\infty}(G_{1}),\del_{1})$ and
$(L^{\infty}(G_{2}),\del_{2})$ with trivial cocycles as follows. Let $\rho: G_{1}\times G_{2}
\to \Omega^{-1}$ be the homeomorphism given by $\rho(g,s)=j(s)i(g)$. Let
$\mathcal{O}=\theta^{-1}(\Omega \cap \Omega^{-1})$ and for $(g,s)\in \mathcal{O}$ define
$\beta_{s}(g)\in G_{2}$ and $\alpha_{g}(s)\in G_{2}$ by
\[\rho^{-1}(\theta(g,s))=(\beta_{s}(g),\alpha_{g}(s)).\] Finally, one can define a
$\ast$-isomorphism $$\tau: L^{\infty}(G_{1})\otimes L^{\infty}(G_{2}) \to L^{\infty}(G_{1})
\otimes L^{\infty}(G_{2})$$ by $\tau(f)(g,s)=f(\beta_{s}(g),\alpha_{g}(s)$. Then, $(\tau,1,1)$
gives a cocycle matching of $(L^{\infty}(G_{1}),\del_{1})$ and $(L^{\infty}(G_{2}),\del_{2})$
with trivial cocycles.

{\sc Example 1.}

\[G=\left\{\begin{pmatrix}
a & b & x\\ c & d & y\\ 0 & 0 & 1
\end{pmatrix}
 | \begin{pmatrix} a & b\\ c & d
\end{pmatrix} \in SL_{2}\R,\ x,y \in \R \right\}
\]
So, $G$ is a Lie-subgroup of $SL_{3}(\R)$.
\[G_{1}=\R^{2},+ \tekst{and} G_{2}=SL_{2}(\R).\]
Further, $i$ maps $G_{1}$ into $G$ in the canonical way and
\begin{displaymath}
i((x,y))=\begin{pmatrix} 1 & 0 & -x\\ -x & 1 & -y+\frac{1}{2}x^2\\ 0 & 0 & 1
\end{pmatrix} \qquad \qquad j(\begin{pmatrix} a & b\\ c & d
\end{pmatrix})=\begin{pmatrix} d & -b &0\\-c&a&0\\0&0&1 \end{pmatrix}
\end{displaymath} Suppose that
\[Q=\begin{pmatrix}
a & b\\ c & d
\end{pmatrix} \in SL_{2}(\R) \; .\]
Then, the mutual actions are given by
\begin{displaymath}
\alpha_{(x,y)}(A)=\begin{pmatrix} a+bx & b\\
c+dx-(a+bx)(ax+b(y+\frac{1}{2}x^2)) & d-b(ax+b(y+\frac{1}{2}x^2))
\end{pmatrix}
\end{displaymath}
and
\begin{displaymath}
\beta_{A}((x,y))=(ax+by+\frac{b}{2}x^2,cx+d(y+\frac{1}{2}x^2)-\frac{1}{2}(ax+b(y+\frac{1}{2}x^2))^{2})
\end{displaymath}
We take trivial cocycles and construct the bicrossed locally
compact quantum group $(M,\del)$. It is not so difficult to show
that $\delta_{1}$ and $\delta_{2}$ are trivial and
$\delta(A,(x,y)) = \mbox{det}\ A = 1$. Therefore, the bicrossed
product is a Kac algebra. One might think that there is a hope to
leave the Kac algebra 'world', if we would work with the general
linear groups (GL) instead of the special linear groups (SL).
Unfortunately, the determinant will be $\alpha$-invariant. So, we
will also find that the bicrossed product is a Kac algebra.

Now, one can construct the infinitesimal Hopf algebra of the bicrossed product quantum group
in the sense of \cite{VV}. It is an algebraic version of the same quantum group.

In this example the infinitesimal Hopf algebra has generators $X$, $Y$, $A$, $B$, $C$ and $D$
satisfying $AD-BC=1$ and the folowing relations
\begin{displaymath}
\begin{array}{lll}
[A,B]=0, & [A,C]=0, & [A,D]=0,\\

[B,C]=0, & [B,D]=0, & [C,D]=0,
\end{array}
\end{displaymath}
\begin{displaymath}
\begin{array}{cc}
\begin{array}{l} [A,X]=B, \\

[B,X]=0,\\

[C,X]=D-A^2,\\

[D,X]=-AB, \end{array}
\begin{array}{l} [A,Y]=0,\\

[B,Y]=0,\\

[C,Y]=-AB,\\

[D,Y]=-B^{2}, \end{array}
\end{array}
\end{displaymath}
\begin{align*}
\del(A)&=A\otimes A + B \otimes C,\\
\del(B)&= B \otimes D + A \otimes B,\\
\del(C)&= C \otimes A + D \otimes C,\\
\del(D)&= D \otimes D + C \otimes B,\\
\del(X)&= 1\otimes X + X \otimes A + Y \otimes C,\\
\del(Y)&= 1 \otimes Y + X \otimes B + Y \otimes D.
\end{align*}

{\sc Example 2.}

Now, we will construct a non-amenable locally compact quantum group that is not a Kac algebra.
\[G_{1}=\{(x,z)|x\in \R, x\neq 0, z\in \C \}\ \tekst{with} (x,z)(y,u)=(xy,z+xu),\]
\[G_{2}=\left\{ \left(
\begin{array}{cc}
a & \bar{c}\\ c & \bar{a}
\end{array}  \right) | a,c\in \C, |a|^2 - |c|^2 = \pm 1
\right\},\]
\[G=\{(2\times 2)-\mbox{matrices over}\ \C\ \mbox{with determinant}
=\pm 1 \} \; .\] Define $Sq(x)=Sgn(x)\sqrt{|x|}$ for all $x \in \R$. Take embeddings $i$ and
$j$ defined by
\begin{displaymath} i:(x,z) \mapsto \frac{1}{Sq(x)} \begin{pmatrix} x & -z\\ 0 &
1  \end{pmatrix} \qquad \qquad j:\begin{pmatrix} a & \bar{c} \\ c & \bar{a}
\end{pmatrix} \mapsto \left( \begin{array}{cc} a & \bar{c} \\ c & \bar{a}
\end{array} \right)^{-1}.
\end{displaymath} The mutual actions are given by
\begin{displaymath}
\alpha_{(x,z)}(a,c)=\frac{D}{Sq(|cz+\bar{a}x|^{2}-|c|^{2})}(\bar{c}\bar{z}+ax,c).
\end{displaymath}
\begin{displaymath}
\beta_{(a,c)}(x,z)=\frac{D}{|x|}(|cz+\bar{a}x|^{2}-|c|^{2},(az+\bar{c}x)(\bar{c}\bar{z}+ax)-
a\bar{c})
\end{displaymath}
with $D=D(x,z,a,c)=\frac{x}{|x|}(|a|^2 - |c|^2)$.

Taking $\cU=\cV=1$, we can construct the bicrossed product locally
compact quantum group $(M,\del)$. Since $\delta$ and $\delta_{2}$
are trivial and $\delta_{1}(x,z)=\frac{1}{x^2}$, we conclude,
using Propositions~2.17 and 4.16  of \cite{VV}, that $(M,\del)$ is
not a Kac algebra, is non-compact and non-discrete. As far as we
know, there was, until now, no example of a non-discrete
non-amenable quantum group that is not a group.

Now, the infinitesimal Hopf $\ast$-algebra is generated as a $\ast$-algebra by normal elements
$A$, $C$ and $Y$, an antiselfadjoint element $X$ and a selfadjoint element $U$ satisfying the
following commutation relations:

\begin{displaymath}
\begin{array}{c}
[A,C]=[A,C^{*}]=0, \qquad  A^{*}A-C^{*}C=U, \qquad U^{2}=1,\\[1ex]
[X,Y]=Y,\\[1ex]
\begin{array}{l@{\qquad}l}
[A,X]=-UACC^*, & [C,X]=-UAA^*C,\\

[A,Y]=2C^{*}-UAA^{*}C^{*}, & [C,Y]=-UA^{*}C^{*}C,\\

[A,Y^{*}]=UA^{2}C, & [C,Y^{*}]=UAC^{2}.
\end{array}
\end{array} \end{displaymath}
Furthermore, the comultiplication is given by
\begin{align*}
\del(A)&=A \otimes A + C^{*}\otimes C,\\
\del(C)&=C \otimes A + A^{*} \otimes C,\\
\del(X)&=X \otimes U(A^{*}A+C^{*}C) + Y \otimes UA^{*}C - Y^{*} \otimes UAC^{*} + 1 \otimes X,\\
\del(Y)&=1 \otimes Y + X \otimes 2UA^{*}C^{*} + Y \otimes U(A^{*})^{2} - Y^{*} \otimes
U(C^{*})^{2}.
\end{align*}

\end{document}